\documentclass[11pt,a4paper]{article}
\usepackage{t1enc}
\usepackage[latin1]{inputenc}
\usepackage[english]{babel}
\usepackage{harvard}
\usepackage{amssymb,stmaryrd,mathabx}

\newcommand{\R}{\mathbb R}

\newcommand{\beq}{\begin{equation}}
\newcommand{\eeq}{\end{equation}}
\newcommand{\beqarr}{\begin{eqnarray}}
\newcommand{\eeqarr}{\end{eqnarray}}
\newcommand{\beqa}{\begin{eqnarray*}}
\newcommand{\eeqa}{\end{eqnarray*}}
\begin{document}
\thispagestyle{empty}

\title{\bf  The changing faces of  the {\em Problem of Space} in the work of Hermann Weyl}
\author{Erhard Scholz\footnote{University of Wuppertal, Faculty of  Math./Natural Sciences, and Interdisciplinary Centre for History and Philosophy of Science, \quad  scholz@math.uni-wuppertal.de}}
\begin{small}
\date{Aug 13, 2016 } 
\maketitle
\begin{abstract}
During his life Weyl approached the {\em problem of space} (PoS)   from various sides. Two aspects stand out as permanent features of his different approaches: the {\em unique  determination of an affine connection} (i.e., without torsion in the  terminology of Cartan) and the question {\em which type of group}  characterizes physical space. The first feature came up  in 1919 (commentaries to Riemann's inaugural lecture)  and played a crucial role in Weyl's work on the PoS in the early 1920s. He defended the central role of affine connections  even in the light of  Cartan's more general framework of connections with torsion. 
In later years, after the rise of the Dirac field, it could have become problematic, but Weyl saw the challenge  posed  to Einstein gravity by  spin coupling primarily in the possibility to allow for  non-metric affine connections. Only after Weyl's death Cartan's approach to  infinitesimal homogeneity and torsion  became revitalized in gravity theories.
\end{abstract}

{ \tableofcontents}
\end{small}

\section*{Introduction}
\addcontentsline{toc}{section}{\protect\numberline{}Introduction}
  According to H. Weyl  three aspects have  to be taken into account for studying  the {\em problem of space} (PoS):  the  extensive medium of the world (``extensives Medium der Aussenwelt''), its metrical structure, and its  content by a material quality  changing from place to place (``materielle Erf\"ullung mit einem von Stelle zu Stelle ver\"anderlichen Quale'') \cite[p. 205]{Weyl:ARP1922b}. The problem has 
to be approached from two sides, by a philosophical investigation and by mathematical analysis. In the early 1920s  Weyl saw  the task of {\em philosophy} in clarifying the distinction and the mutual relationships between the  three aspects mentioned above.   In addition to this {\em mathematics}  had to ``search for  correct knowledge of the {\em essence of space}  and of the {\em spatial structure}'' as far as   it  can be ``described  quantitatively, in logico-arithmetical relations''.\footnote{``F\"ur den Mathematiker handelt es sich darum, das quantitativ, in logisch-arithmetischen Relationen Erfa\ss{}bare am Wesen des Raumes und der r\"aumlichen Struktur richtig zu erkennen und  mit den Hilfsmitteln der Logik, Arithmetik und Analysis auf seine einfachsten Gr\"unde zur\"uckzuf\"uhren'' \cite[p. 206]{Weyl:ARP1922b}. }  Note Weyl's double description of {\em essence} and {\em structure};  he considered them as  complementary aspects of the  concept of space.

Such a characterization was given by him in the phase 1921 to 1923 when Weyl developed his program of the {\em mathematical analysis of the problem of space} in a well defined,  sense. In the following we shall denote it by {\em PoS$_{21-23}$}. 
It dealt  with the  question of how to generalize the  Helmholtz-Lie  analysis of the homogeneity conditions of classical space to the new context of relativistic physics.  Weyl insisted on the necessity to reformulate homogeneity in terms of differentiable manifolds endowed with linear groups operating in the infinitesimal neighbourhoods (in modernized  language, operating on the tangent spaces), or between them. He gave very general conceptually motivated conditions and analyzed their consequences.  His result was an infinitesimal group structure  typical for the automorphisms of his generalized differential geometry, ``pure infinitesimal geometry'', developed  in 1918 ({\em Weylian metric}).  The contribution of Weyl to the problem of space has found much attention in the history and philosophy of mathematics.
Here it  will be dealt with  from a specific point of view only; for more aspects and finer  details the reader may consult the literature.\footnote{Among many \cite{Scheibe:Spacetime,Sigurdsson:Diss,Coleman/Korte:DMV,Bernard:Idealisme,Bernard:Barcelona,Scholz:Weyl_PoS,%
Scholz:Weyl/Cartan}.  Weyl's study of the PoS$_{21-23}$ was the guiding axis of the conference of which this book arose; see in particular the contributions to this volume by A. Roca-Rosell and  C. Lobo. \label{fn lit PoS}} 

This was not the only situation in which Weyl addressed the problem of space. In a more general sense  this problem  was a recurrent theme in his thought all over his life. Weyl  hit upon  the PoS  (in the wider sense) in different contexts and looked  at it from different angles.  
The present paper puts  Weyl's discussion of the PoS in a wider perspective (section \ref{section faces}), but it would be far beyond its scope to deal with all these  different facets in some detail. 

Here we shall concentrate on selected topics which show how Weyl used  context dependent   relative  a priori elements  which he considered constitutive for determining the structure   of space, or even for grasping its ``essence''. Two conceptual features stand out among  them: (i)  the characterization of homogeneity by means of group structures  and (ii) the core role assigned by Weyl to a uniquely determined affine (torsion free) connection for admissible space structures. Both features appear prominently in Weyl's PoS$_{21-23}$ but also, in different form, in other encounters of him with the problem of space. 
We shall discuss how Weyl dealt with the problem of homogeneity after the rise of general relativity (section \ref{section homogeneity})  and contrast it with Cartan's answer to the question  (section \ref{section Cartan's homogeneity}).  That could have given reasons for Weyl to revise  the central role of uniquely determined affine connections as a kind of relative a priori for physical geometry, but it did not (section \ref{section affine connection}).  The next challenge  was posed by Dirac's spinor fields in  general relativity. In the light of later developments (Einstein-Cartan theory of gravity) one might expect that it could have become  a problem for Weyl's affine connection principle already in the 1930/40s. But that was not the case.  Why, will be shortly discussed in section \ref{section spinor fields}, before we come to a final evaluation (section \ref{section discussion}). 

Of course Weyl, like the other mathematicians of the 20th century, often used the terminology of ``space'' at other places in a wider sense than above, sometimes in a purely mathematical context. But in the framework of this paper the notion of {\em space} is nearly always used in the more restricted sense of a mathematical space structure which serves, or at least may serve, as a candidate for grasping  physical space or space-time, the ``extensive medium of the external world'' in Weyl's word. For the abbreviation PoS this is always the case.

\section{\small The multiple faces  of the PoS \label{section faces}}
As already mentioned we can give here only  a short survey of  different contexts and different forms in which Weyl met the problem of space. The following list of topics and contexts  may serve  as an orientation (main publications indicated in brackets):

\begin{enumerate}
\item  Modernized presentation of  the classical problem of space in the sense of Helmholtz and Lie in the first chapter of {\em Raum - Zeit - Materie} \cite{Weyl:RZM5}
	\item Specification of  Riemannian metrics (``Pythagorean nature'' of metric) in the wider class of Finsler metrics in Riemann's approach to geometry \cite{Riemann/Weyl:Hypothesen}
		\item Mathematical analysis of the problem of space, PoS$_{21-23}$, \cite{Weyl:ARP1923}
		\item Characterization of $\R^3$ or $S^3$ (the three-dimensional sphere) by combinatorial invariants.  Topological space forms and their characterization by discrete groups (operating on the universal covering space) \cite[pp. 16ff.]{Weyl:RiemannsIdeen}
			\item Introduction of differential structure on continuous manifolds, in particular with regard to its justification fo the concept of physical space  \cite[p. 12]{Weyl:RiemannsIdeen}
		\item  Cartan's general concept of infinitesimal geometries   \cite[pp. 38ff.]{Weyl:RiemannsIdeen}, \cite{Weyl:1929[82],Nabonnand:Cartan/Weyl}
				\item Similarities and congruences as an exemplary case for the distinction of mathematical automorphisms and physical automorphisms of space \cite{Weyl:PMNEnglish,Weyl:Hs91a:31}
			\item  Specific role of Lorentz/Poincar\'e group and the dimension 4 of spacetime \cite{Weyl:Hs91a:31}\footnote{See the contribution of S. de Bianchi, this volume.}
		\item Finally Weyl's  considerations on the possible role of  non-metric affine connections for the dynamics of  spinor fields in general relativity \cite{Weyl:1950}
\end{enumerate}

Under the items 1., 3., 6., 7.  Weyl dealt with the question of how to characterize the  homogeneity of the respective spatial structures. Here different versions of automorphism groups played a prominent role. In his discussion of the items 2., 3. and 9. Weyl's conviction that a proper space structure in the sense of PoS carries a uniquely determined affine (i.e. torsion free linear) connection stood out. In 2., 3. the postulate of a uniquely determined affine connection was not questioned at all; in 9. he subjected it to a check whether it could be defended in the light of relativistic spinor fields. Item 7. and 8. have been  discussed elsewhere.\footnote{\cite{Scholz:Weyls_search}}
 In the following sections this  topic will be dealt in more detail. The topological space problem (item 4) and the problem of differentiability (item 5) have only been brought up at  isolated occasions by Weyl; they cannot be discussed in this paper.



\section{\small Homogeneity characterizations  given by Weyl  \label{section homogeneity}}
In the first three editions of {\em Raum - Zeit - Materie}\,  the classical space problem of the 19th century,  posed and answered by Helmholtz, Lie, Engel (Weyl added Hilbert, {\em Grundlagen der Geometrie}, app. IV) was mentioned by Weyl only in passing \cite[pp. 86, 264, 1st to 3rd ed.]{Weyl:RZM}. Even in 
the fourth edition published in 1921, in which he already included a first sketch of his own thoughts on the PoS$_{21-23}$, Weyl did not go into more details \cite[pp. 86, 289]{Weyl:RZM4}. Only {\em after} having finished his own analysis he gave a more extended presentation of the classical solution in the fifth edition \cite[p. 100]{Weyl:RZM5}. He did not  use the  concepts of ``rigid body'' and ``free mobility'', which had become problematic with the advent of relativity theory, but expressed Helmholtz' postulates  abstractly in terms of   group theoretical constraints for the {\em homogeneity} of classical space. He  rephrased Helmholtz's axioms of free mobility by conditions which are now called {\em simple flag transitivity} of the homogeneity group.\footnote{``Man kommt so zu der folgenden Formulierung des Homogeneit\"atspostulats
im $n$-dimensionalen Raum: Es soIl m\"oglich sein, mit
Hilfe einer zur Gruppe {\bf G} geh\"origen Abbildung ein System $\Sigma$ inzidenter
Richtungselemente der $0^{ten}$ bis $(n-1)^{ten}$ Stufe in ein gleichartiges, beliebig
vorgegebenes System $\Sigma'$ \"uberzuf\"uhren; aber die Identit\"at soll unter
den Abbildungen von {\bf G} die einzige sein, welche ein derartiges System $\Sigma$
inzidenter Richtungselemente festl\"a\ss{}t'' \cite[p. 100]{Weyl:RZM5}. } Similar in \cite[5th lecture]{Weyl:ARP1923}. 

 Weyl's presentation stripped Helmholtz' analysis from the latter's intention of founding his conditions on  supposedly {\em factual} conditions (``Thatsachen'') lying at the basis of any  empirical measurement.  He did not claim to give a historically precise account of Helmholtz' thoughts,  in fact his  passage  read as though  Helmholtz had started from an investigation of the a priori  conditions of the  homogeneity of physical space.\footnote{``Von einem tieferen, gruppentheoretischen Gesichtspunkt aus hat Helmholtz
zuerst die Homogeneitatsfrage gestellt. Helmholtz setzt nicht die
G\"ultigkeit des Pythagoreischen Lehrsatzes im Unendlichkleinen, ja nicht
einmal die Me\ss{}barkeit der Linienelemente voraus; er spricht allein von
dem wahren Grundbegriff der Geometrie, der Gruppe {\bf G} der kongruenten
Abbildungen des Raumes.'' \cite[p. 100]{Weyl:RZM5}}
Weyl's reconstruction of Helmholtz's and Lie/Engels' approach to the space problem was   systematic, not historical. It assimilated  the classical PoS to a perspective which prepared the way to a type of analysis which Weyl pursued in his own program between  1921 and 1923. From such a perspective  the classical PoS  seemed to have been   the question of an a priori characterization  of the homogeneity of space. It led to an answer  which allowed to introduce an invariant definite quadratic differential form of constant curvature and  thus to the classical spaces of Euclidean and non-Euclidean geometry.\footnote{``Es ist eine wunderbare gruppentheoretische Tatsache, die von Helmholtz, strenger und allgemeiner von
S. Lie bewiesen wurde, da\ss{}  die einzigen dieser Bedingung gen\"ugenden
Gruppen {\bf G} die Gruppen {\bf G}$_{\lambda}$ \ldots  [sind]'' ibid.  By {\bf G}$_{\lambda}$  Weyl denoted the congruence groups of hyperbolic, parabolic, or elliptic geometry (in the terminology introduced by F. Klein).}

Such a type of a priori characterization could no longer  be apodictic, like  Kant's  a priori judgements had been (or, at least, had been claimed to be). It  no longer consisted of necessary judgements, but  rather of well founded postulates,  if possible of the {\em best founded} ones, which characterize the conceptually {\em possible}  in a specified context. In this function it  still has an a priori character {\em in distinction from} empirical determinations but   only {\em relative to the latter and to the theoretical context}.
 With a change of context and/or more refined empirical knowledge,  formerly well established a priori conditions can  become obsolete and may have to be  revised, in agreement with the analysis given by \citeasnoun[59ff.]{Friedman:Reconsidering}.  

Such a revision was the goal of Weyl's analysis of the problem of space (1921--1923). By several reasons simple flag transitivity could no longer be upheld as a feature characterizing the homogeneity of space. Firstly special relativity had integrated space proper and time to spacetime as the ``extensive medium of the world''. That  
destroyed flag transitivity because it now became necessary to account for the qualitative difference of  timelike and spacelike directions. Moreover general relativity, Einstein's theory of gravity,  broke with  the paradigm of constant curvature and  made curvature dependent on  the distribution of matter and energy,  thus giving it an  a posteriori character.  Therefore Weyl, and a little later Cartan,  posed the question of homogeneity of spacetime  anew, in  forms adapted to the context of general relativity. Both came to different, only partially overlapping answers which we shall discuss in the following.

Weyl started from a conceptual analysis of what he considered the most general, minimal conditions for congruence geometry founded on infinitesimal structures like those he had proposed in his purely infinitesimal geometry of 1918. In his investigations 1921--1923 he wanted to dig deeper and to motivate, or even derive, a  generalized  metrical structure from congruence and similarity concepts in the infinitesimal and a generalized homogeneity principle. In a move which he presented as a conceptual analysis of the idea of congruence in the infinitesimal he established conditions for (linear) groups characterizing {\em congruence} by generalized ``rotations'' in the  infinitesimal neighbourhoods of each point of  spacetime. An intuitive idea of {\em homogeneity} then demanded that the type of group (more precisely the conjugation class in the general linear group) was equal for any two points. But Weyl did not use the  terminology ``homogeneous/homogeneity'' in this context; he rather spoke of the fixed {\em nature} of the group (the conjugation class) which could be expressed by pointwise changing ``orientations'' of the group (members of the class)\cite[p. 48]{Weyl:ARP1923}.\footnote{Similarly in the 5th edition of {\em Raum - Zeit - Materie}, where he characterized the ``nature'' of a Riemannian  metric by its signature and its ``orientation'' by the point dependent value of the respective quadratic differential form \cite[p. 102]{Weyl:RZM5}. In the 4th edition (translated into English and French)  he still fought more indirectly  with the problem that space as  ``a form of phenomena \ldots is necessarily homogeneous'', while the Riemannian metric is not  \cite[pp. 96ff.]{Weyl:RZMengl}.}

Moreover a transition between neighbouring points had to be specified by a {\em linear connection} which is compatible with the {\em similarities} with regard to the ``rotations'' (more technically with the normalizer of the rotations). All this was justified by  by what Weyl considered an a priori  {\em analysis of the  concepts} involved \cite[p. 49]{Weyl:ARP1923}. 

But his was not all; Weyl added:
\begin{quote}
I now come to the {\em synthetic part} in the {\sc Kant}ian sense. The task is now to    formulate precisely the postulate,  up to now only  indicated, which finally determines  the type of rotation group which is characteristic for the real world. (ibid.)\footnote{``Ich komme jetzt zum {\em synthetischen} 
Teil im {\sc Kant}ischen Sinne. Da gilt es, das fr\"uher angedeutete Postulat
pr\"azis zu formulieren, das die f\"ur die wirkliche Welt charakteristische Art
der Drehungsgruppe festlegen soll'' \cite[p. 49]{Weyl:ARP1923}, emph. in original.}
\end{quote}
The  {\em synthetic} component of his a priori justification of infinitesimal congruence structures consisted of a two-part postulate, the first of which demanded a kind of wide adaptability to a posteriori distribution structure of matter (postulate of ``free deformability''), the second one was the postulate of {\em unique determination of a compatible affine connection}. Later the first part  turned out to be mathematically redundant,\footnote{\cite{Scheibe:Diss} } leaving the second part as the mathematical and philosophical core of Weyl's synthetic a priori of the PoS$_{21-23}$.\footnote{For more details see the literature cited in fn. \ref{fn lit PoS}.}

Mathematically it was  crucial for constraining the groups which could serve as candidates for infinitesimal congruences so strongly that in the end  Weyl could show that  only the generalized orthogonal groups (of arbitrary signature) satisfy the constraints ({\em main theorem of Weyl's PoS}). 
A philosophical analogy of this principle with  intersubjectivity relations in practical philosophy   was expressed and emphasized  by Weyl in his Barcelona lectures \cite[p. 46]{Weyl:ARP1923}. The nature of this   analogy 
is being analyzed by N. Sieroka (this volume)   and related to Weyl's exchange with F. Medicus and his reading of Fichte. 

The existence of a uniquely determined affine connection remained a stable feature in Weyl's understanding of  a {\em good} geometric structure designed for representing {\em space} from 1919 onward, although the mathematical feasibility of  it might have became doubtful after Cartan's answer to the homogeneity challenge of general relativity and even stronger after  the advent of relativistic spinor fields. Only much later, in  the years between 1948 and 1950,  Weyl  subjected it to an investigation of its empirical acceptability. His   result was that this part of his a priori may be sustained even {\em relative} to a context including general relativistic spinor fields (see section \ref{section spinor fields}).

With regard to  homogeneity in the modern (general relativistic) context  Weyl performed  a tight-rope walk. Clearly, space ``as a form of phenomena \ldots is necessarily homogeneous'' \cite[p. 96]{Weyl:RZMengl}, but Riemannian spaces are not. 
He solved, or circumvented,  the problem by arguing that (simply connected) neighbourhoods of any two points are diffeomorphic and the ``nature'' of the metric remains the same for all points.\footnote{Speaking in global terms, Weyl would surely  have assumed a transitive diffeomorphism group of the spacetime manifold.}  Both belong to the a priori determinations of space, while the metrical structure is fixed by a posteriori, empirically given factual conditions (or contingent model assumption, we might add).\footnote{Similar in  \cite[p. 87]{Weyl:PMNEnglish}.}  
For him the manifold and the nature of the metric, both invariant under diffeomorphisms, expressed the generalized idea of homogeneity in the relativistic context. But he did not speak  of ``generalized homogeneity'', he rather reserved the terminology of homogeneity for the classical situation of  metrically homogeneous spaces \cite{Weyl:RZM,Weyl:ARP1923,Weyl:PMNEnglish}. 

In spite of this terminological decision,  the view  that the 
diffeomorphisms of the  spacetime manifold are part of the {\em physical automorphisms} of general relativistic field theories  persisted in Weyl's thought.  At the time of  preparing the different editions of {\em Raum - Zeit - Materie} (1918--1923) Weyl argued for such an understanding by means of the plasticine analogy discussed in J. Bernard's contribution, this volume. This was an intuitive, ``didactical'', way for expressing the more general postulate that under a physical automorphism the field structures are ``dragged along'' with the diffeomorphisms which can thus be considered as dynamical symmetries in present physicists' terminology. Otherwise they would not be able to preserve the (a posteriori) field structures and the metric induced by them. Weyl insisted on such a generalized understanding of homogeneity in a talk  on physical and mathematical automorphisms given in the late 1940s \cite{Weyl:Hs91a:31}.\footnote{Cf. \cite{Scholz:Weyls_search}}

\section{\small  Cartan's concept of infinitesimal homogeneity  \label{section Cartan's homogeneity}}
Elie Cartan approached the problem of homogeneity in a different way.
 In 1922 he  presented a new type of infinitesimal geometry to the public {\em Sur une g\'en\'eralisation de la notion de courbure \ldots} \cite{Cartan:1922[58]}. He had developed the necessary tools  (differential forms and Lie group theory)  over a long time  and elaborated the basic ideas for his new geometry during the year 1921 in an interplay of  differential geometry, Einstein's gravity theory, and  the brothers Cosserat's  generalized theory of elasticity.\footnote{\cite{Cogliati:Stockholm,Nabonnand:2016Cartan,Scholz:Weyl/Cartan};  for a modern mathematical introduction to Cartan geometry see \cite{Sharpe:Cartan_Spaces}.}  In the years to come he would expand his approach to a broad program of generalized infinitesimal geometries later called {\em Cartan spaces}.  In his survey talk at the International Congress of Mathematicians in Toronto, 1924, Cartan
motivated the approach by indicating that  general relativity was confronted with 
\begin{quote}
  \ldots the paradoxical task of
interpreting in a non-homogeneous universe \ldots the multiple experiences made by observers who believed  in the homogeneity of this universe \cite[p. 86]{Cartan:1924Toronto}.
\end{quote}
 Although general relativity  had helped to induce a first step towards  bridging  the gap between (non-homogeneous) Riemannian geometry and Euclidean geometry (homogeneous in the sense of F. Klein)  by motivating T. Levi-Civita's   concept of infinitesimal parallelism, he did not yet see the gap  closed.\footnote{``Or, c'est le d\'eveloppement m\^eme de la th\'eorie de la relativit\'e,
li\'ee par l'obligation paradoxale d'interpr\'eter dans et par un
Univers non homog\`ene les r\'esultats de nombreuses exp\'eriences
faites par des observateurs qui croyaient \`a l'homogen\'eit\'e de cet
Univers, qui permit de combler en partie le foss\'e qui s\'eparait
les espaces de Riemann de l'espace euclidien. Le premier pas
dans cette voie fut l'oeuvre de M. Levi-Civita, par
l'introduction de sa notion de parall\'elisme.'' \cite[p. 86]{Cartan:1924Toronto})}

Cartan alluded to the Kleinean understanding of homogeneity and indicated the idea underlying his approach:
\begin{quote} . . .while a Riemannian space does not possess absolute homogeneity,
it possesses a kind of infinitesimal homogeneity; in the
immediate neighbourhood of a point it can be assimilated to a
Euclidean space. (ibid.)\footnote{``\ldots  si un espace de Riemann ne poss\`ede pas une homog\'en\'eit\'e absolue, il poss\`ede cependant une sorte d'homog\'en\'eit\'e infinit\'esimale; au voisinage immediat d'un point donn\'e il est assimilable \`a une espace euclidien'' \cite[p. 85]{Cartan:1924Toronto}.}
\end{quote}
That is, Cartan wanted to implement  {\em infinitesimal homogeneity} in his new generalized spaces, in addition to their infinitesimal (generalized) rotational symmetries.  Weyl, as we have seen, translated homogeneity in the new, relativistic context to the possibility of comparing the neighbourhoods of any two points of spacetime (not infinitesimally close ones)  with each other, which  boilt down to considering  the diffeomorphisms of the underlying manifold as part of the automorphisms of the  geometric structure.

Cartan's program thus consisted in an  infinitesimalization of the Kleinean concept of geometry as a homogeneous space $S$  in the sense of $S\cong G/H$ with a main group $G$ and generalized rotations (isotropy group) $H \subset G$.\footnote{In the Euclidean case $G\cong \R^n \rtimes SO(n,\R)$ with $H=SO(n,\R)$,  $S=G/H\cong \R^n$.} 
 Cartan considered the infinitesimal version of the groups, i.e. the  corresponding Lie algebras $\mathfrak g = Lie \, G, \mathfrak h = Lie \, H$ and $\mathfrak g = \mathfrak l \oplus \mathfrak h$, and assimilated the quotient of the two  with the infinitesimal neighbourhoods in the manifold $M$, such that $   T_p M \cong \mathfrak g / \mathfrak h \cong \mathfrak{l}$,  for any point $p$ of $M$. His crucial symbolical tools   were ensembles of differential forms, which can be read as differential forms with values in the respective Lie algebras. With their help he introduced a generalized type of connection, now called {\em Cartan connection}, which led to two kinds of curvature effects.\footnote{Cf. \cite{Sharpe:Cartan_Spaces}.} 

The curvature with respect to the generalized rotations $\mathfrak h$ corresponded to the Riemann curvature (in  slightly different form) which was well known at the time, while the curvature with respect to the  generalized translations $\mathfrak l$ was a new effect. Cartan  called  it {\em torsion} because in the context of the generalized elasticity theory in the sense of the Cosserats it could be    interpreted    as a rotational  momentum in the medium. 
If translated to the coordinate notation of differential geometry (the calculus of Ricci and Levi-Civita) Cartan's connection could, in many cases, be expressed in the form of a linear connection $\Gamma$ with coefficients $\Gamma_{\mu \nu}^{\lambda}$ which are no longer symmetric in the lower indices. In fact the condition $\Gamma_{\mu \nu}^{\lambda}\neq \Gamma_{\nu \mu}^{\lambda}$ is equivalent to { non-vanishing torsion}.\footnote{$\Gamma_{\mu \nu}^{\lambda}-\Gamma_{\nu \mu}^{\lambda} = T_{\nu \mu}^{\lambda}$ is the  {\em torsion} tensor, i.e., the translational curvature expressed in coordinate coefficients.  }

Cartan did not try to argue on a philosophical level as explicitly as Weyl did, but  his conceptual analysis of the ``paradoxical task'' posed by general relativity may be read as  establishing a new type  of  a priori framework for physical geometry, different from Weyl's in   PoS$_{21-23}$.
Cartan's new relativ a priori was wider than Weyl's in two respects. His framework allowed a larger variety of infinitesimal isotropy and homogeneity groups than Weyl's.
Moreover, in his view it would not  appear  natural  to consider the  existence and uniqueness of a compatible affine connection (torsion zero) as a ``synthetic'' a priori of the physical space concept.  Cartan  reformulated Weyl's PoS in his framework, but  he dealt with it  from a mathematical point of view rather then of a philosophical one, and with a slightly different outcome.\footnote{Cf. \cite{Scholz:Weyl/Cartan}.}

 Weyl  responded to Cartan's proposal of a large class of infinitesimal geometries, but at the beginning he was not convinced that the wider  perspective was helpful for extending the a priori framework of physical geometry. In the correspondence between him and Cartan he expressed doubts even with regard to the specific geometrical achievements of Cartan's generalization, although at the end both authors came to a basic acceptance of the other's viewpoint \cite{Nabonnand:Cartan/Weyl}.\footnote{For a survey  see \cite{Scholz:Weyl/Cartan}; a more refined evaluation of the correspondence is being prepared by C. Eckes and P. Nabonnand.} 
In his contribution to the Lobachevsky anniversary volume 
 written in 1925, but published only posthumously, Weyl discussed Cartan's approach and acknowledged that it allowed a ``far-reaching generalization of infinitesimal geometry''  \cite[p. 38]{Weyl:RiemannsIdeen};  in particular:
\begin{quote} \ldots it achieves the natural  widest possible range of concept formation which still allows to establish a theory of curvature in analogy to Riemann's.\footnote{``Und darauf beruht wohl \"uberhaupt die mathematische Bedeutung seines allgemeinen Schemas: es erreicht den nat\"urlichen
weitesten Umfang der Begriffsbildung, welche die Aufstellung einer
Kr\"ummungstheorie analog der Riemannschen noch erm\"oglicht'' \cite[p. 39]{Weyl:RiemannsIdeen}. }  
\end{quote}
Thus Weyl acknowledged Cartan's generalization of the curvature concept, but without mentioning  that it carries the potential to undermine the central role of the affine connection, which  he considered as the most important  part of the ``synthetic'' a priori of the space concept.\footnote{In the same article he reiterated that he still stood to the content of his PoS$_{21-23}$: ``Das neue gruppentheoretische Raumproblem, das vom
Standpunkt der Relativit\"atstheorie an Stelle des Helmholtz-Lie'schen
tritt, glaube ich in meiner Schrift "Mathematische Analyse
des Raumproblems" (1923, Vorlesung 7 und 8) formuliert und gel\"ost
zu haben.'' \cite[p. 37]{Weyl:RiemannsIdeen}}

\section{\small  Affine connection, synthetic a priori or just a special condition among others?   \label{section affine connection}}
Shortly after Levi-Civita's invention of the infinitesimal parallel displacement in Riemannian geometry Weyl introduced affine connections  as a concept  of its own for differential geometry, which allowed to talk about {\em infinitesimal parallel displacements} in allusion to Levi-Civita's terminology but independent of  the structure of a Riemannian metric and without reference to an embedding into a higher dimensional Euclidean space \cite{Weyl:InfGeo}. He simply demanded that (1) for any vector $\xi$ attached to a point $p$ of the manifold the change  induced by parallel displacement  from  $p$ to an infinitesimal close point $p'$ 
depends linearly on the vector  $\widearrow{pp'}$,
and (2) if for two points $p_1,p_2$, both infinitesimally close to $p$, the parallel displacement of $\widearrow{pp_1}$ along $\widearrow{pp_2}$ leads to the  end point $p_{21}$ and the displacement of $\widearrow{pp_2}$ along $\widearrow{pp_1}$  to  $p_{12}$, then $p_{12}= p_{21}$: ``The result is an infinitesimally small parallelogram.''\footnote{``Es entsteht eine unendlich kleine Parallelogrammfigur'' \cite[p. 7]{Weyl:InfGeo}.} 

 In the paper directed to mathematicians  Weyl argued   conceptually, sometimes even in a   philosophical style. He  presented the task of geometry being ``to fathom out  the essence of the metrical concepts''.\footnote{Die Geometrie ``ergr\"undet, was im Wesen der metrischen Begriffe liegt'' \cite[p. 2]{Weyl:InfGeo}. In the paper presenting his purely infinitesimal geometry to physicists (as the geometrical background for his  unified field theory)  Weyl introduced the affine connection in more concrete form and axiomatically \cite[p. 32]{Weyl:GuE}, \cite[p. 26]{Weyl:GuE_English}.}
He thus understood the conditions (1) and (2) as  postulates which arise from an {\em analysis of the concept} of infinitesimal parallel displacement. In the philosophical  language used   in  \cite{Weyl:ARP1923} the generalized affine connection  resulted from   an a  priori  conceptual analysis.

 A little later, in his commentaries to Riemann's inaugural lecture \cite{Riemann/Weyl:Hypothesen}, he pondered on the question of how the Riemannian metric could be specified among the wider class of Finsler metrics (which had a striking  a priori justification in being homogeneous with regard to rescaling). He conjectured that the Riemannian metrics are just those which are compatible with a uniquely determined affine connection.\footnote{``Bei der fundamentalen Bedeutung, die nach den neueren
Untersuchungen (\ldots)
dem affinen Grundbegriff der infinitesimalen Parallelverschiebung
eines Vektors f\"ur den Aufbau der Geometrie
zukommt, erhebt sich insbesondere die Frage, ob die Mannigfaltigkeiten
der Pythagoreischen Raumklasse die einzigen
sind, welche die Aufstellung dieses Begriffs erm\"oglichen und
welche dementsprechend nicht blo\ss{} eine Metrik, sondern
auch affinen Zusammenhang besitzen. Die Antwort lautet
wahrscheinlich bejahend, ein Beweis daf\"ur ist aber bisher
nicht erbracht worden.'' \cite[p. 27]{Riemann/Weyl:Hypothesen}.}
D. Laugwitz would later call this conjecture  {\em  Weyl's first problem of space} -- and answered it positively \cite{Laugwitz:ARP}.

Weyl was convinced of the  fundamental import of the  principle of a uniquely determined affine connection already in 1919; in the following we shall speak about it as  Weyl's {\em affine connection principle}.  Its crucial role for deriving the 
main theorem of the PoS$_{21-23}$   made it  so convincing for Weyl that in 1922/23 he even  raised it to the status of a {\em  ``synthetic'' a priori}. 

Considered from a wider perspective this was neither self-evident nor imperative (in distinction to a Kantian understanding of synthetic a priori). Only a few weeks after  Weyl's Barcelona lectures in February 1922,  E. Cartan presented his first public note on  generalized spaces to the Paris {\em Acad\'emie des sciences}  \cite{Cartan:1922[58]}. Of  course Weyl  could not know about it at the time of his lectures, nor apparently  while preparing them for publication, but in hindsight it could have became clear to him  that  Cartan's generalized spaces also opened the pathway towards  a different relative a priori  for relativistic spacetime.

 For Cartan  the difference was not so much of a philosophical nature, but mathematically it was clear to him from  the outset. 
If  a parallel displacement,  $\Gamma(\xi, \widearrow{pp'}) = \xi'-\xi$, is expressed by  connection coefficents $\Gamma^{\lambda}_{\mu\nu}$ with regard to a coordinate basis of the infinitesimal neighbourhoods (the tangent spaces), the closing condition (2) boils down to the symmetry of the coefficients, $\Gamma^{\lambda}_{\mu\nu}=\Gamma^{\lambda}_{\nu\mu}$.    Cartan's torsion tensor expressed in coordinate coefficents, on the other hand, bcomes $T^{\lambda}_{\mu\nu}= \Gamma^{\lambda}_{\nu\mu}-\Gamma^{\lambda}_{\mu\nu}$; Weyl's  condition (2) thus amounts  to vanishing torsion.\footnote{Cartan discussed this point in  his investigation of Weyl's PoS  in slightly different terms \cite[\S 3]{Cartan:1923PoS}.}

It was not easy for Weyl and Cartan to disentangle their different points of view on their differences with regard to  generalized spaces, although  Cartan could treat the mathematical aspects of  Weyl's PoS$_{21-23}$  quite easily as a special case of his methods and he acknowledged Weyl's deep philosophical analysis, but without discussing it from his side \cite{Cartan:1923PoS}.  Weyl, on the other hand, found it difficult to grasp the subtleties of Cartan's approach,  while he soon understood   the wider mathematical generality of the  latter's approach and acknowledged it (see above).  In 1929  he even adapted certain aspects of Cartan's  approach  for his proposal to of a general relativistic version of  Dirac's electron theory \cite{Weyl:1929Dirac}. In the same year he  gave a survey talk on Cartan's theory to the Princeton group of differential geometry and mathematical physics. There he argued that for a proper geometric usage one has to add certain restrictions to Cartan's scheme, among them the exclusion of  torsion \cite[p. 210]{Weyl:1929[82]}. 

 Cartan had developed and extended his methodology in the course of the 1920s in very general terms, with infinitesimal Klein spaces of many different types and  even without assuming that, pointwise, the  infinitesimal quotient group (Liealgebra) $\mathfrak l \cong \mathfrak g/\mathfrak h$ can be identified with the tangent spaces of the manifold (not even the dimension needed to be the same). 
For  Weyl  it seemed indispensable  that for a proper geometrical usage of Cartan's general scheme one had  to  assume pointwise identifications of   $\mathfrak l $, which he called the ``tangent plane'' denoted by ``$T_P$'' (sic!),   with the infinitesimal neighbourhoods of the manifold, the tangent spaces in the later terminology.\footnote{In the modern understanding of  Cartan spaces this is an indispensable property  inbuilt in the definition of a {\em Cartan gauge},  e.g. \cite[p. 174]{Sharpe:Cartan_Spaces}. Note that Weyl's ``$T_P$'' stood for for $\mathfrak l $ in the function of  what would now be considered the tangent space of the translative subgroup in the fibre direction.}
 He called this an ``embedding'' of ``$T_P$''  into the manifold  and added additional restrictions motivated by the  specific geometrical structure considered. He spoke of ``special manifolds'',   in particular with regard to projective and conformal structures.\footnote{Weyl  drew upon  the results of the Princeton school of differential geometry, A.L. Eisenhart, O. Veblen, T.Y Thomas, intending to build bridges between the French (E. Cartan) and the US (Princeton school) traditions in differential geometry.}
The specialization conditions contained, in particular, the ``invariantive restriction to require that our manifold \ldots is without torsion'' \cite[p. 210]{Weyl:1929[82]}. 

In the  following correspondence Cartan insisted that his research program did not need such restrictions and deplored that it was not fairly represented by Weyl's survey. The ensuing exchange of letters  centered on the role of  ``embedding'' of $\mathfrak l$ (Weyl's ``$T_P$) and the specialization conditions. Although torsion played only a subordinate role, the correspondence shows once again that Cartan considered  torsion zero  as a technical specialization condition among others without particular conceptual import \cite{Nabonnand:Cartan/Weyl}.

\section{\small  The challenge of spinor fields in the 1930-40s  \label{section spinor fields}}
In the light of later developments in Einstein-Cartan theory of gravity (see final section) it ought to be added that, to my knowledge, in the 1930s   neither  Cartan nor Weyl considered the question whether the general relativistic  Dirac electron field with non-vanishing spin (spin = $\frac{1}{2}$) might undermine Weyl's affine principle  from the physical, perhaps even from an empirical side. Even  in the late 1940s, when Weyl started to analyze the consequences of an independent spin coupling of the Dirac field to gravity, he did not pose the question whether the affine principle had to be given up in favour of Cartan's view (torsion $\neq 0$). 
 He rather  chose an approach which Einstein had studied in the middle of the 1920s, in which an affine connection and the metric of a generalized Lagrangian were varied independently  \cite{Einstein:1925}. He thus  {\em relaxed the condition of metricity} of the connection  (calling it a ``mixed'' theory) rather than that of vanishing torsion. 
 For a Lagrangian of the electron field with a Dirac term, a spin term and a mass term   Weyl found that ``by the influence of matter a slight discordance between affine connection and metric is created'' \cite[p. 288 and equs. (2), (3)]{Weyl:1950}. 

This interesting observation clearly  posed a challenge to Weyl's affine connection principle.  But  ``by somewhat laborious calculations'' (not presented in the paper) he was able to show  that by adding a small term of the form $-12\pi G\,  l_2$  to the  Lagrangian ($l_2$ a quartic scalar invariant in the 4-component spinor field $\psi$) the metric theory became equivalent to the   mixed theory (without the additional term).\footnote{$l_2= (\overline{\psi}_3\psi_1+\overline{\psi}_4\psi_2)(\overline{\psi}_1\psi_3+\overline{\psi}_2\psi_4)$.  The ``laborious calculations'' seem to have been presented in the manuscript (not preserved) for  the publication \cite{Weyl:Coupling[140]}.  K. Chandrasekharan, the editor of Weyl's {\em Gesammelte Abhandlungen} remarks about this paper  that ``due to typographical errors, it is incomprehensible'' \cite[vol.  4, p. 285, footnote]{Weyl:GA}. It  was therefore not reprinted in the edition. }
 He   concluded:
\begin{quote} To this extent then there is a a complete equivalence between the mixed metric-affine  and the purely metric conception of gravity. \cite[p. 288]{Weyl:1950} 
\end{quote}
The metric theory of gravity could be upheld by a small addition to the Lagrangian even in the  light of an  electron field's spin coupling to gravity. In this way the challenge posed by Dirac spinor fields to Weyl's  affine connection principle was  neutralized and the relative a priori of the uniquely determined metric affine connection successfully defended.

\section{\small Late endorsement for Cartan's infinitesimal homogeneity principle and general discussion  \label{section discussion}}
In his PoS$_{21-23}$ Weyl tried to found a new conceptual framework for space and time that lived up to the challenge of the theories of relativity, like the  homogeneous spaces of the late 19th century had done with regard to classical physics. The  challenge arose from physical theories in which Einstein had evaluated both, empirical and theoretical knowledge in a quite specific sense. Regarding physical concepts Einstein was an ingenious innovator (role of time, space, simultaneity, equivalence principle, metric as gravitational field etc.), but with regard to the mathematical theories, he  had built with conceptual material inherited from contemporaneous mathematics (Riemannian geometry, Ricci-Levi Civita's absolute calculus).  Weyl intended  to go beyond the constraints of the inherited  and to ``fathom out'' (as he said in the above quote) the minimal  ingredients of  spacetime concepts,  necessary for obtaining  infinitesimal congruence structures which were able to build a bridge between the general notion of a differentiable manifold and specific metrical determinations. The latter ought to be able to adapt as flexible as possible to contingent  distributions of matter and fields. He did so in what he considered an a priori move of  conceptual analysis,  as he said in open allusion to the   Kantian terminology, and found that he had to add the ``synthetic'' affine connection principle. Although he tried to motivate the synthetic principle by analogy with considerations of practical philosophy (respect for a ``common good'') the finally decisive motivation of the principle lay in its success for deriving the main theorem of the PoS, not different from what one would ordinarily expect from  the axioms  of a mathematical theory. 

Weyl was well aware that his   a priori was different from Kant's; in particular it was no longer apodictic and made sense only in the  context of relativistic physics and the  horizon of new differential geometric structures on manifolds. It thus  was relative with regard to the theoretical context and open  for potential revisions, like the classical understanding of homogeneity had been. In Weyl's view this did not make the striving for a well understood a priori obsolete. In his view the role of a priori statements had changed  from being necessary judgements to {\em well motivated  possible} concepts and structures. But their function with regard to more specific theoretical and empirical knowledge  remained.  In his view ``physics projects what is given onto the background of the possible'' \cite[p. 220]{Weyl:PMNEnglish} and mathematics  explores the conceptually possible.\footnote{``The dual nature of reality accounts for the fact that we cannot design a theoretical image of being except upon the background of the possible. Thus the four-dimensional continuum of space and time is the field of the a priori existing possibilities of coincidences.'' \cite[p. 231]{Weyl:PMNEnglish}} 
 
Substituting Kant's  necessary a priori judgements by  the ``background of the possible'' points  also towards another shift:  the relative a priori need no longer be uniquely determined.   In our case study we have come across a possible loss of uniqueness,  the {\em underdetermination of the relative a priori},  by comparing  Weyl's analysis of the PoS with Cartan's generalized spaces. The latter's conceptual motivation was the implementation of  infinitesimal homogeneity in addition to the global  homogeneity achieved by structure dragging diffeomorphisms. That  gave them the potential  for becoming a competing a priori structure for relativistic space concepts. In our discussion of Weyl's different approaches to the space problem this potentiality did not materialize, apparently because Cartan did not like to argue too much on the philosophical level and Weyl did not see  imperative reasons to give up his affine connection principle. It needed a change of generations and  deep conceptual as well as technical studies of physicists  before the a priori  potential of Cartan spaces became apparent. 
A detailed account of this next shift would need a publication of its own (or  more). Here we have to content ourselves with an outline.\footnote{For a rich collection of sources with detailed commentaries from the theoretical physics side see \cite{Blagojevic/Hehl}.}
 
The success of conceiving electromagnetism as a $U(1)$-gauge field induced physicists to study field theoretic consequences of point dependent infinitesimal symmetries of other groups. In the terminology of physics the symmetries were ``localized''.\footnote{See A. Afriat's paper, this volume.}  
Most striking and best known is the case of at first strong, later weak isospin $SU(2)$ (Yang/Mills) and the later generalizations to gauge field theories in elementary particle physics. But also  the Poincar\'e group, the  symmetry group of special relativity, was  ``localized''  independently and nearly simultaneously by T. Kibble and D. Sciama \cite{Sciama:1962,Kibble:1961}. From the point of view of physics  the conserved currents of infinitesimal symmetries supplied by the Noether theorems played a crucial heuristic role for this  research. For the Poincar\'e group $\R^4 \rtimes SO(3,1)$ that led to considering the spin current, the Noether current with regard to the Lorentz rotations, as an additional source for the gravitational field, supplementing    energy-momentum, the  current with regard to the translation group.\footnote{To be more precise: Sciama  presupposed an Einsteinean background and gained spin as an additional current,  modifying Einstein gravity to what was later called Einstein-Cartan gravity.  Kibble, on the other hand, started from localizing the symmetries of Minkowski space and considered different Lagrangians, the simplest of which led to Einstein-Cartan theory \cite[p. 106]{Blagojevic/Hehl}.} 

The physical idea of ``localizing'' the translations of the Poincar\'e group together with the Lorentz rotations was very close to Cartan's idea to implement infinitesimal homogeneity in addition to infinitesimal isotropy in his concept of generalized spaces, although  neither Kibble nor Sciama noticed the kinship of their studies with Cartan's geometrical framework.  This was brought into the open  by the work of F. Hehl and  A. Trautman.\footnote{\cite{Hehl:Habil,Trautman:1973,Hehl_ea:1976} and others.} 
Then it also became clear that the simplest Lagrangian in Kibble's approach, as well as in  Sciama's,  is equivalent to the one discussed  by Cartan  in passing, when he showed what his approach could contribute to understand and to extend Einstein's theory of gravity \cite[\S83]{Cartan:1923[66]}. It is now being called {\em Einstein-Cartan gravity} (EC).\footnote{Cf. \cite{Trautman:2006,Hehl:Dennis}.} 

  Einstein-Cartan gravity modifies Einstein's general relativity only to a tiny degree; for  vanishing spin it reduces to the latter. Moreover,  outside of spinning matter field the torsion is zero and the influence of spin on the metric can be taken into account by a small modification of the energy-momentum source of the Einstein equation,\footnote{Cf. \cite[p. 406]{Hehl_ea:1976}, \cite[p. 194]{Trautman:2006}.}
 similar to what Weyl had found for the non-metricity induced by spin in the ``mixed'' theory.
Only for mass densities more than $10^{38}$ times the one  of a neutron star, respectively a  nucleon mass compressed to $10^6$ Planck lengths,  which signifies energy densities at the hypothetical grand unificaton scale of elementary particle interactions, the experts expect EC ``to overtake'' Einstein's general relativity.\footnote{\cite[p. 194]{Trautman:2006} \cite[p. 108]{Blagojevic/Hehl}.}

These seemingly technical results of modified gravity are important in our context, because they show that Cartan's {\em principle of infinitesimal homogeneity} has  finally become important in foundational studies of gravity. In the last third of the 20th century it has  turned into a {\em relative apriori} for relativistic spacetime theories with, at least, the same right as Weyl's affine connection principle and alternative to the latter. It  even would have the advantage of being  closer to what Weyl called the ``physical automorphisms'' of modern physics 
by fitting well to the Noether current paradigm for infinitesimal symmetries, prominent in contemporary physics \cite{Weyl:Hs91a:31}.\footnote{Cf. \cite{Scholz:Weyls_search}.}

On the other hand, there is (still?) no actual empirical evidence which would force us to revise Weyl's analysis of the PoS$_{21-23}$ and to relegate his  affine connection principle from the status of a relative apriori to an empirical principle, valid only in ``weak'' field constellations. In this sense, we seem to be here in the situation of a presumably {\em temporal underdetermination} of the relative a priori principles of Weyl and Cartan.\footnote{The so-called teleparallel version of  Cartan geometric gravity rearranges the coordination between Noether currents and dynamical equations: Energy-momentum becomes the source of translational curvature and the spin current of the rotational curvature, rather than the other way round as in EC gravity. In oral communications D. Lehmkuhl has indicated that teleparallel gravity may be an interesting example of a possibly principled (in contrast to temporal) underdetermination of empirically equivalent gravity theories.}
This seems to be another feature of present a priori  structures, at least as important as their  being established in the context of wider scientific results and being open to revision with them.

\small
 \bibliographystyle{apsr}
  \bibliography{a_lit_hist,a_lit_mathsci,a_lit_scholz}

\end{document}